\documentclass{article}
\usepackage{amsmath}
\usepackage{amsmath,amssymb,latexsym,color}
\usepackage[mathscr]{eucal}
\newtheorem{thm}{Theorem}[section]
\newtheorem{prop}[thm]{Proposition}
\newtheorem{lemma}[thm]{Lemma}

\newtheorem{example}{Example}[section]
\newtheorem{defin}{Definition}[section]

\newtheorem{remark}{Remark}[section]

\newcommand{\proof}{{\it Proof.\quad}}
\newcommand{\qed}{\hfill\Box\medskip}

\usepackage{CJK}
\begin{document}
\begin{CJK*}{GBK}{song}

\renewcommand{\baselinestretch}{1.15}
\title{Orthogonal graphs over Galois rings of odd characteristic}

\author{
Fenggao Li$^{\rm a}$\quad Jun Guo$^{\rm b}$ \quad Kaishun
Wang$^{\rm c}$\footnote{Corresponding author.\newline  Email
addresses: fenggaol@163.com (F. Li), guojun$_-$lf@163.com (J. Guo), wangks@bnu.edu.cn (K. Wang)}  \\
{\footnotesize  $^{\rm a}$ \em  College of Mathematics, Hunan Institute of Science and Technology, Yueyang, 414006, China}\\
\footnotesize $^{\rm b}$ \em  Math. and Inf. College, Langfang Teachers' College, Langfang, 065000, China\\
\footnotesize $^{\rm c}$ \em Sch. Math. Sci. {\rm \&} Lab. Math. Com. Sys., Beijing Normal University, Beijing, 100875,  China}
\date{}
\maketitle

\begin{abstract}
Assume that $\nu$ is a positive integer and $\delta=0, 1$ or $2$. In this paper we introduce the orthogonal graph $\Gamma^{2\nu+\delta}$ over a Galois ring of odd characteristic and prove that it is arc transitive.
Moreover, we compute its parameters as a quasi-strongly regular graph. In particular, we show that $\Gamma^{2+\delta}$ is a strongly regular graph and $\Gamma^{2\nu+1}$ is a strictly Deza graph when $\nu\geq 2$.

\medskip


\end{abstract}

\section{Introduction}

Let $p$ be an odd prime number, and $s,m$ be positive integers. We use $R$ to denote the Galois
ring of characteristic $p^{s}$ and cardinality $p^{sm}$. Note that $R$ is a finite field with $p^{m}$ elements if $s=1$, and $R$ is the ring of residue classes of $\mathbb{Z}$ modulo its ideal $p^{s}\mathbb{Z}$ if $m=1$.
By \cite{Wan2}, there exists a unit $\xi$ of multiplicative order $p^{m}-1$ such that every element $a$ of $R$ can be written uniquely as
$a=a_{0}+a_{1}p+\cdots+a_{s-1}p^{s-1}$,
where each $a_{i}$ belongs to $\{0,1,\xi,\ldots,\xi^{p^{m}-2}\}$. Note that $a$ is a unit if and only if $a_{0}\ne 0$. It follows that the unit group $R^*$ of $R$ has the size $(p^{m}-1)p^{m(s-1)}$.

Suppose $\nu$ is a positive integer and $\delta=0,1$ or 2. Let $R^{2\nu+\delta}$ be the set of $(2\nu+\delta)$-tuples $(a_{1},a_{2},\ldots,a_{2\nu+\delta})$ of elements in $R$ such that $a_{i}\in R^{\ast}$ for some $i$.
Denote by $e_i $ the element of $R^{2\nu+\delta}$ whose $i$th coordinate is 1 and all other coordinates are 0.
Define $(a_1,a_2,\ldots,a_{2\nu+\delta})\sim (b_1,b_2,\ldots,b_{2\nu+\delta})$
if there exists a $\lambda\in R^{\ast}$ such that $(a_1,a_2,\ldots,a_{2\nu+\delta})=\lambda (b_1,b_2,\ldots,b_{2\nu+\delta})$.
Then $\sim$ is an equivalence relation on $R^{2\nu+\delta}$.
Write $[a_1,a_2,\ldots,a_{2\nu+\delta}]$ for the equivalence class containing $(a_1,a_2,\ldots,a_{2\nu+\delta})$,
and let $\tilde{R}^{2\nu+\delta}$ be the set of all equivalence classes.
By \cite{Cao1}, any $(2\nu+\delta)\times (2\nu+\delta)$ inverse symmetric matrix over $R$ is cogredient to
$$
\mathrm{S}_{2\nu+\delta,\Delta}=\left(\begin{array}{ccc}0&I^{(\nu)}&\\[-0.1cm]
                                      I^{(\nu)}&0&\\[-0.1cm]
                                      &&\Delta\end{array}\!\right),
$$
where
$$
\Delta=\left\{\begin{array}{ll}\phi\,\mbox{(disappear)} &\mbox{if }\delta=0,\\[0.1cm]
                  (1)\mbox{ or } (z) &\mbox{if }\delta=1,\\[0.1cm]
                  \mathrm{diag}(1,-z) &\mbox{if }\delta=2,
\end{array}\right.
$$
$I^{(\nu)}$ is the identity matrix of order $\nu$, and $z$ is a fixed non-square element of $R^{\ast}$.
The {\em orthogonal graph} over $R$, denoted by $\Gamma^{2\nu+\delta}$, is the graph whose vertex set is
$\{[\alpha]\in \tilde{R}^{2\nu+\delta}\mid\alpha \mathrm{S}_{2\nu+\delta,\Delta}\alpha^{t}=0\,\}$
and two vertices $[\alpha]$ and $[\beta]$ are adjacent if $\alpha\mathrm{S}_{2\nu+\delta,\Delta}\beta^{t}\in R^{\ast}$.
The {\em orthogonal group} of degree $2\nu+\delta$ over $R$ with respect to $\mathrm{S}_{2\nu+\delta,\Delta}$, denoted by $\mathrm{O}_{2\nu+\delta}(R)$, consists of all $(2\nu+\delta)\times (2\nu+\delta)$ matrices $T$ over $R$ satisfying $T\mathrm{S}_{2\nu+\delta,\Delta}T^t=\mathrm{S}_{2\nu+\delta,\Delta}$.
Note that $\mathrm{O}_{2\nu+\delta}(R)$ is an automorphism group of $\Gamma^{2\nu+\delta}$, and there is an action of $\mathrm{O}_{2\nu+\delta}(R)$ on $\Gamma^{2\nu+\delta}$ defined by
$$
\begin{array}{ccc}
           V(\Gamma^{2\nu+\delta})\times\mathrm{O}_{2\nu+\delta}(R)& \longrightarrow &V(\Gamma^{2\nu+\delta})\\[0.1cm]
            ([x_1,x_2,\ldots,x_{2\nu+\delta}],T) & \longmapsto &[(x_1,x_2,\ldots,x_{2\nu+\delta})T],
\end{array}
$$
where $V(\Gamma^{2\nu+\delta})$ is the vertex set of $\Gamma^{2\nu+\delta}$.

The collinearity graph of the classical polar space is a family of well known strongly regular graphs (cf. \cite{BL, Gu, TW, Wan1}).
In the symplectic case, the complement of this graph was generalized to some special local rings in \cite{LWG, MP, MP2}. In this paper we shall study the orthogonal graph over a Galois ring of odd characteristic with $s\geq 2$.

In Section 2 we show that $\Gamma^{2\nu+\delta}$ is arc transitive, and list all representatives of suborbits of $\mathrm{O}_{2\nu+\delta}(R)$ on $\Gamma^{2\nu+\delta}$. In Section 3 we compute all parameters of $\Gamma^{2\nu+\delta}$ as a quasi-strongly regular graph, and show that $\Gamma^{2+\delta}$ is a strongly regular graph and $\Gamma^{2\nu+1}$ is a strictly Deza graph for $\nu\geq 2$.

\section{Suborbits}

In this section we prove that $\Gamma^{2\nu+\delta}$ is arc transitive, and list all representatives of suborbits of $O_{2\nu+\delta}(R)$ on $\Gamma^{2\nu+\delta}$.

Let $E_{i}$ be the $\nu\times 1$ matrix whose $(i,1)$-entry is 1 and all other entries are 0, and $E_{ij}$ be the $\nu\times 2$ matrix whose $(i,j)$-entry is 1 and all other entries are 0. We write $\tilde{\alpha}_{\delta}=\emptyset (\mathrm{disappear})$, $a_{2\nu+1}e_{2\nu+1}$ or $a_{2\nu+1}e_{2\nu+1}+a_{2\nu+2}e_{2\nu+2}$ according to $\delta=0,1$ or 2, respectively. Write
$$
H_{2\nu+\delta,\Delta}(P)=\left(\begin{array}{ccc}
             I^{(\nu)}&-\frac{1}{2}P\Delta P^{t}&-P\\
             0&I^{(\nu)}&0\\
             0&\Delta P^{t}&I^{(\delta)}\\
            \end{array}
       \right),\;\;\;\;
T_{2\nu+\delta,\Delta}(P)=\left(\begin{array}{ccc}
             I^{(\nu)}&0&0\\
             -\frac{1}{2}P\Delta P^{t}&I^{(\nu)}&-P\\
             \Delta P^{t}&0&I^{(\delta)}\\
            \end{array}
       \right),
$$
$$
K_{2\nu+\delta}(Q)=\left(\begin{array}{ccc}
             I^{(\nu)}&Q&\\
             &I^{(\nu)}&\\
             &&I^{(\delta)}\\
            \end{array}
       \right),\;\;
Q_{t}(x_{1},\ldots,x_{t-1})=\left(\begin{array}{cccc}
             0&-x_{1}&\ldots&-x_{t-1}\\
             x_{1}&&&\\
             \vdots&&\\
             x_{t-1}&&&\\
            \end{array}
       \right),
$$
where $P$ is a $\nu\times\delta$ matrix and $Q$ is a $\nu\times\nu$ alternate matrix over $R$.
Then $H_{2\nu+\delta,\Delta}(P)$, $T_{2\nu+\delta,\Delta}(P)$ and $K_{2\nu+\delta}(Q)$ belong to $\mathrm{O}_{2\nu+\delta}(R)$.

Let $S$ be an $m\times m$ inverse symmetric matrix over the finite field $\mathbb{F}_{q}$. If for any vector $x\in\mathbb{F}_{q}^{m}$ such that $xSx^{t}=0$ we have necessarily $x=0$, then $S$ is said to be {\em definite}.

\begin{lemma}\label{lem2.1}
For any vertex $[a_1,a_2,\ldots,a_{2\nu+\delta}]$ of $\Gamma^{2\nu+\delta}$, there exists an $i\in\{1,2,\ldots,2\nu\}$ such that $a_i\in R^{\ast}$.
\end{lemma}
\proof
If $\delta=0,$ the result is directed. If $\delta=1$, then $2(a_{1}a_{\nu+1}+\cdots+a_{\nu}a_{2\nu})+\Delta a_{2\nu+1}^{2}=0$ implies the desired result. Now suppose $\delta=2$. If $a_i\notin R^{\ast}$ for each $i\in\{1,2,\ldots,2\nu\}$, then $2(a_{1}a_{\nu+1}+\cdots+a_{\nu}a_{2\nu})+a_{2\nu+1}^{2}-za_{2\nu+2}^{2}=0$ implies $\bar{a}_{2\nu+1}^{2}-\bar{z}\bar{a}_{2\nu+2}^{2}=\bar{0}$ in the quotient field $R/pR$.
By \cite[Chapter 6]{Wan}, the matrix $\mathrm{diag}(\bar{1},-\bar{z})$ over $R/pR$ is definite. So $\bar{a}_{2\nu+1}=\bar{a}_{2\nu+2}=\bar{0}$, i.e., $a_{2\nu+1},a_{2\nu+2}\notin R^{\ast}$, which are impossible. Hence, the desired result follows.
$\qed$

\begin{thm}\label{thm3.1}
The orthogonal graph $\Gamma^{2\nu+\delta}$ is arc transitive.
\end{thm}
\proof
First we prove that $\Gamma^{2\nu+\delta}$ is vertex transitive. Let $[\alpha]=[a_1,a_2,\ldots,a_{2\nu+\delta}]$ be any vertex of $\Gamma^{2\nu+\delta}$. It suffices to show that there exists an element of $\mathrm{O}_{2\nu+\delta}(R)$ carrying $[\alpha]$ to $[e_{1}]$.
By Lemma~\ref{lem2.1} there are the following two cases to be considered.

{\em Case\,1}: $a_{i}\in R^{\ast}$ for some $i\in\{1,2,\ldots,\nu\}$. Note that there exists a $\nu\times\nu$ inverse matrix $S_{11}$ over $R$ such that $(a_{1},\ldots,a_{\nu})S_{11}=(1,0,\ldots,0)$. Then
$S_{1}=\mathrm{diag}(S_{11},(S_{11}^{-1})^{t},I^{(\delta)})\in \mathrm{O}_{2\nu+\delta}(R)$ carries $[\alpha]$ to
$[e_{1}+b_{\nu+1}e_{\nu+1}+\cdots+b_{2\nu}e_{2\nu}+\tilde{\alpha}_{\delta}]$, where $(b_{\nu+1},\ldots,b_{2\nu})=(a_{\nu+1},\ldots,a_{2\nu})(S_{11}^{-1})^{t}$.
If $\delta=0$, then $K_{2\nu+\delta}(Q_{\nu}(b_{\nu+2},\ldots,b_{2\nu}))$ carries $[\alpha S_{1}]$ to
$[e_{1}]$.
If $\delta=1$, then $H_{2\nu+1,\Delta}(a_{2\nu+1}E_{1})$ carries $[\alpha S_{1}K_{2\nu+1}(Q_{\nu}(b_{\nu+2},\ldots,b_{2\nu}))]$ to $[e_{1}]$. If $\delta=2$, then $H_{2\nu+2,\Delta}(a_{2\nu+1}E_{11}+a_{2\nu+2}E_{12})$
carries $[\alpha S_{1}K_{2\nu+2}(Q_{\nu}(b_{\nu+2},\ldots,b_{2\nu}))]$ to $[e_{1}]$.

{\em Case\,2}: $a_{1},\ldots,a_{\nu}\notin R^{\ast}$, $a_{\nu+i}\in R^{\ast}$ for some $i\in\{1,2,\ldots,\nu\}$. Then
$$
\left(\begin{array}{ccc}
             0&I^{(\nu)}&\\
             I^{(\nu)}&0&\\
             &&I^{(\delta)}\\
            \end{array}
       \right)
\in\mathrm{O}_{2\nu+\delta}(R)
$$
carries $[\alpha]$ to $[a_{\nu+1}e_{\nu+1}+\cdots+a_{2\nu}e_{2\nu}+a_{1}e_{1}+\cdots+a_{\nu}e_{\nu}+\tilde{\alpha}_{\delta}]$,
reduced to {\em Case\,1}.

Now we prove that $\Gamma^{2\nu+\delta}$ is arc transitive.
If $[\alpha]$ is adjacent to $[e_{1}]$, then $[\alpha]$ is of the form
$[\alpha]=[a_{1},\ldots,a_{\nu},1,a_{\nu+2},\ldots,a_{2\nu+\delta}]$. It suffices to show that there exists an element of
the stabilizer $G_{[e_1]}$ of $[e_1]$ in $O_{2\nu+\delta}(R)$ carrying $[\alpha]$ to $[e_{\nu+1}]$.
If $\delta=0$, then
$$
\left(
    \begin{array}{cccc|cccc}
     1&&&&&&&\\
     a_{\nu+2}&1&&&&&&\\
     \vdots&&\ddots&&&&&\\
     a_{2\nu}&&&1&&&&\\ \hline
     a_{1}&-a_{2}&\ldots&-a_{\nu}&1&-a_{\nu+2}&\ldots&-a_{2\nu}\\
     a_{2}&&&&&1&&\\
     \vdots&&&&&&\ddots&\\
     a_{\nu}&&&&&&&1\\
    \end{array}
  \right)
$$
is a desired element. Similar to the case of $\delta=0$, the desired result follows for the remaining cases.
$\qed$

In order to give all representatives of suborbits of  $\mathrm{O}_{2\nu+\delta}(R)$ on $\Gamma^{2\nu+\delta}$, by Theorem~\ref{thm3.1} it suffices to consider the action of $G_{[e_{1}]}$ on $\Gamma^{2\nu+\delta}$.

\begin{lemma}~\label{lem3.2}
Let $[\alpha]=[a_{1},a_{2},\ldots,a_{2\nu+\delta}]$ be a vertex of $\Gamma^{2\nu+\delta}$, where $a_{i}\in R^{\ast}$ for some $i\in\{2,3,\ldots,2\nu\}$. Then there exists an element of $G_{[e_{1}]}$ carrying $[\alpha]$ to
$$
\begin{array}{ll}
[e_{2}] & \mbox{if $\nu=1$},\\[0.1cm]
[e_{\nu+1}],[e_{2}]\;\; \mbox{or}\;\;[e_{2}+p^{r}e_{\nu+1}] & \mbox{if $\nu\geq 2$},\\
\end{array}
$$
where $0<r<s$.
\end{lemma}
\proof
If $a_{\nu+1}\in R^{\ast}$, then $[e_{1}]$ is adjacent to $[\alpha]$. By Theorem~\ref{thm3.1}, there exists an element of $G_{[e_{1}]}$ carrying $[\alpha]$ to $[e_{\nu+1}]$. Now assume that $a_{\nu+1}\notin R^{\ast}$. Then $\nu\geq 2$.
By Lemma~\ref{lem2.1} there are the following two cases to be considered.

{\em Case\,1}: $a_{i}\in R^{\ast}$ for some $i\in\{2,\ldots,\nu\}$. Then there exists a $(\nu-1)\times(\nu-1)$ inverse matrix $T_{11}$ over $R$ such that $(a_{2},\ldots,a_{\nu})T_{11}=(1,0,\ldots,0)$. Write $(b_{\nu+2},\ldots,b_{2\nu})=(a_{\nu+2},\ldots,a_{2\nu})(T_{11}^{-1})^{t}$ and
$$
T_{1}=\mathrm{diag}(1,T_{11},1,(T_{11}^{-1})^{t},I^{(\delta)})
\left(
  \begin{array}{ccccc}
   1&&&&\\
   &I^{(\nu-1)}&&Q_{\nu-1}(b_{\nu+3},\ldots,b_{2\nu})&\\
   &&1&&\\
   &&&I^{(\nu-1)}&\\
   &&&&I^{(\delta)}\\
  \end{array}
\right).
$$
Then $T_{1}\in G_{[e_{1}]}$ and $[\alpha T_{1}]=[a_{1}e_{1}+e_{2}+a_{\nu+1}e_{\nu+1}+b_{\nu+2}e_{\nu+2}+\tilde{\alpha}_{\delta}]$.
Write $c_{\nu+2}=b_{\nu+2}+a_{1}a_{\nu+1}$. Then
$$
T_{2}=\left(\begin{array}{ccc|ccc|c}
      1&&&&&&\\
      -a_{1}&1&&&&&\\
      &&I^{(\nu-2)}&&&&\\ \hline
      &&&1&a_{1}&&\\
      &&&&1&&\\
      &&&&&I^{(\nu-2)}&\\ \hline
      &&&&&&I^{(\delta)}\\
\end{array}\right)\in G_{[e_{1}]}
$$
carries $[\alpha T_{1}]$ to $[e_{2}+a_{\nu+1}e_{\nu+1}]$, $[e_{2}+a_{\nu+1}e_{\nu+1}+c_{\nu+2}e_{\nu+2}+\tilde{\alpha}_{1}]$ or $[e_{2}+a_{\nu+1}e_{\nu+1}+c_{\nu+2}e_{\nu+2}+\tilde{\alpha}_{2}]$ according to $\delta=0,1$ or 2, respectively.
When $a_{\nu+1}\ne 0$ we may write $a_{\nu+1}=xp^{r}$ for some $x\in R^{\ast}$ and $0<r<s$.
If $\delta=0$, then $[\alpha T_{1}T_{2}]=[e_{2}]$ when $a_{\nu+1}=0$; and $\mathrm{diag}(x,I^{(\nu-1)},x^{-1},I^{(\nu-1)})\in G_{[e_{1}]}$ carries $[\alpha T_{1}T_{2}]$ to $[e_{2}+p^{r}e_{\nu+1}]$ when $a_{\nu+1}=xp^{r}$.
Note that $H_{2\nu+1,\Delta}(a_{2\nu+1}E_{2})\in G_{[e_{1}]}$
carries $[\alpha T_{1}T_{2}]$ to $[e_{2}+a_{\nu+1}e_{\nu+1}]$ when $\delta=1$, and $H_{2\nu+2,\Delta}(a_{2\nu+1}E_{21}+a_{2\nu+2}E_{22})\in G_{[e_{1}]}$ carries $[\alpha T_{1}T_{2}]$ to $[e_{2}+a_{\nu+1}e_{\nu+1}]$ when $\delta=2$, reduced to the case similar to $\delta=0$.

{\em Case\,2:} $a_{2},\ldots,a_{\nu}\notin R^{\ast}$, $a_{\nu+i}\in R^{\ast}$ for some $i\in\{2,\ldots,\nu\}$. Then
$$
\left(\begin{array}{cc|cc|c}
      1&&&&\\
      &0&&I^{(\nu-1)}&\\ \hline
      &&1&&\\
      &I^{(\nu-1)}&&0&\\ \hline
      &&&&I^{(\delta)}\\
\end{array}\right)\in G_{[e_{1}]}
$$
carries $[\alpha]$ to $[a_{1}e_{1}+a_{\nu+2}e_{2}+\cdots+a_{2\nu}e_{\nu}+a_{\nu+1}e_{\nu+1}+a_{2}e_{\nu+2}+\cdots+a_{\nu}e_{2\nu}+\tilde{\alpha}_{\delta}]$,
reduced to {\em Case\,1}.
$\qed$

\begin{lemma}~\label{lem3.3}
Let $[\alpha]=[a_{1},a_{2},\ldots,a_{2\nu+\delta}]$ be a vertex of $\Gamma^{2\nu+\delta}$ distinct to $[e_{1}]$,
where $a_{i}\notin R^{\ast}$ for each $i\in\{2,\ldots,2\nu\}$.
If $\nu=1$, then there exists an element of $G_{[e_{1}]}$ carrying $[\alpha]$ to
$$
\begin{array}{ll}
[1,-\frac{1}{2}\Delta p^{2r},p^{r}] & \mbox{if $\delta=1$},\\[0.1cm]
[1,-\frac{1}{2}p^{2r},p^{r},0], [1,\frac{1}{2}zp^{2r},0,p^{r}]\;\;\mbox{or}\;\; [1,-\frac{1}{2}(b^{2}-z)p^{2r},bp^{r},p^{r}] & \mbox{if $\delta=2$},\\
\end{array}
$$
where $b\in R^{\ast}$ and $0<r<s$.
If $\nu\geq 2$, then there exists an element of $G_{[e_{1}]}$ carrying $[\alpha]$ to
$[e_{1}+p^{r}e_{2}]$ or $[e_{1}+p^{r}e_{2}-ap^{r+t}e_{\nu+1}+ap^{t}e_{\nu+2}]$, where $0<r\leq t<s$ and $a\in\{1,z\}$.
\end{lemma}
\proof
Note that $(\nu,\delta)\ne (1,0)$ and $a_{1}\in R^{\ast}$. We may assume that $a_{1}=1$.

First consider the case of $\nu=1$. If $\delta=1$, then $[\alpha]=[1,-\frac{1}{2}\Delta x^{2}p^{2r},xp^{r}]$ for some $x\in R^{\ast}$ and $0<r<s$. Then $\mathrm{diag}(x,x^{-1},1)\in G_{[e_{1}]}$ carries $[\alpha]$ to $[1,-\frac{1}{2}\Delta p^{2r},p^{r}]$. Now assume $\delta=2$. Then $[\alpha]=[1,-\frac{1}{2}x^{2}p^{2r},xp^{r},0]$, $[1,\frac{1}{2}zx^{2}p^{2r},0,xp^{r}]$ or $[1,-\frac{1}{2}(x^{2}p^{2r}-zy^{2}p^{2t}),xp^{r},yp^{t}]$, where $x,y\in R^{\ast}$ and $0<r,t<s$.
If $[\alpha]=[1,-\frac{1}{2}x^{2}p^{2r},xp^{r},0]$ or $[1,\frac{1}{2}zx^{2}p^{2r},0,xp^{r}]$, then $\mathrm{diag}(x,x^{-1},I^{(2)})\in G_{[e_{1}]}$ carries $[\alpha]$ to $[1,-\frac{1}{2}p^{2r},p^{r},0]$ or $[1,\frac{1}{2}zp^{2r},0,p^{r}]$ in respective case.
If $[\alpha]=[1,-\frac{1}{2}(x^{2}p^{2r}-zy^{2}p^{2t}),xp^{r},yp^{t}]$, we can assume $r=t$. In fact, if $r\ne t$, by
\cite[Lemma~1.28]{Wan} we may choose $c,d\in R^{\ast}$ with $c^{2}-d^{2}z=1$. Let
$$
T=\left(
    \begin{array}{cc}
      c&d\\
      dz&c\\
    \end{array}
  \right).
$$
Then $\mathrm{diag}(I^{(2)},T)\in G_{[e_{1}]}$ carries $[\alpha]$ to
$$
\begin{array}{ll}
[1,-\frac{1}{2}(x^{2}p^{2r}-y^{2}zp^{2t}),(cx+dyzp^{t-r})p^{r},(dx+cyp^{t-r})p^{r}] & \mbox{if $r<t$},\\[0.1cm]
[1,-\frac{1}{2}(x^{2}p^{2r}-y^{2}zp^{2t}),(dyz+cxp^{r-t})p^{t},(cy+dxp^{r-t})p^{t}] & \mbox{if $r>t$},\\
\end{array}
$$
where $cx+dyzp^{t-r},dx+cyp^{t-r}\in R^{\ast}$ if $r<t$, and $dyz+cxp^{r-t},cy+dxp^{r-t}\in R^{\ast}$ if $r>t$, as desired.
When $r=t$, $\mathrm{diag}(y,y^{-1},I^{(2)})\in G_{[e_{1}]}$ carries $[\alpha]$ to $[1,-\frac{1}{2}(b^{2}-z)p^{2r},bp^{r},p^{r}]$,
where $b=xy^{-1}\in R^{\ast}$.

Now consider the case of $\nu\geq 2$. We may assume that $(a_{2},\ldots,a_{\nu},a_{\nu+2},\ldots,a_{2\nu})\ne (0,\ldots,0)$. In fact, if
$a_{i}=0$ for some $i\in\{2,\ldots,2\nu\}\setminus\{\nu+1\}$, then $\delta=1$ or 2.
When $\delta=1$ we have that $a_{2\nu+1}\ne 0$ and $T_{2\nu+1,\Delta}(E_{2})\in G_{[e_{1}]}$
carries $[\alpha]$ to $[e_{1}+\Delta a_{2\nu+1}e_{2}+a_{\nu+1}e_{\nu+1}+a_{2\nu+1}e_{2\nu+1}]$ with $\Delta a_{2\nu+1}\ne 0$, as desired.
When $\delta=2$ we have $a_{2\nu+1}\ne 0$ or $a_{2\nu+2}\ne 0$. If $a_{2\nu+1}\ne 0$, then $T_{2\nu+2,\Delta}(E_{21})\in G_{[e_{1}]}$
carries $[\alpha]$ to $[e_{1}+a_{2\nu+1}e_{2}+a_{\nu+1}e_{\nu+1}+a_{2\nu+1}e_{2\nu+1}+a_{2\nu+2}e_{2\nu+2}]$, as desired.
If $a_{2\nu+2}\ne 0$, then $T_{2\nu+2,\Delta}(E_{22})\in G_{[e_{1}]}$
carries  $[\alpha]$ to $[e_{1}-za_{2\nu+2}e_{2}+a_{\nu+1}e_{\nu+1}+a_{2\nu+1}e_{2\nu+1}+a_{2\nu+2}e_{2\nu+2}]$ with $-za_{2\nu+2}\ne 0$, as desired.
Now suppose $(a_{2},\ldots,a_{\nu},a_{\nu+2},\ldots,a_{2\nu})\ne (0,\ldots,0)$. Without loss of generality we may assume that
$a_{2}R+\cdots+a_{\nu}R=p^{r}R$ and $a_{\nu+2}R+\cdots+a_{2\nu}R=p^{r_{1}}R$ with $0<r\leq r_{1}\leq s$.
Then there exists a $(\nu-1)\times (\nu-1)$ inverse matrix $T_{11}$ over $R$ such that
$(a_{2},\ldots,a_{\nu})T_{11}=(p^{r},0,\ldots,0)$.
Write $(a_{\nu+2},\ldots,a_{2\nu})(T_{11}^{-1})^{t}=(b_{\nu+2},\ldots,b_{2\nu})$. Then $T_{1}=\mathrm{diag}(1,T_{11},1,(T_{11}^{-1})^{t},I^{(\delta)})\in G_{[e_{1}]}$ and
$$
[\alpha T_{1}]=[e_{1}+p^{r}e_{2}+a_{\nu+1}e_{\nu+1}+b_{\nu+2}e_{\nu+2}+\cdots+b_{2\nu}e_{2\nu}+\tilde{\alpha}_{\delta}].
$$
When $(b_{\nu+3},\ldots,b_{2\nu})\ne (0,\ldots,0)$ we may assume that $b_{\nu+3}R+\cdots+b_{2\nu}R=p^{r_{2}}R$ for $0<r_{2}<s$. Then $r_{1}\leq r_{2}$ and there exists a $(\nu-2)\times (\nu-2)$ inverse matrix $T_{22}$ over $R$ such that
$(b_{\nu+3},\ldots,b_{2\nu})T_{22}=(p^{r_{2}},0,\ldots,0)$.
Write $T_{2}=\mathrm{diag}(I^{(2)},(T_{22}^{-1})^{t},I^{(2)},T_{22},I^{(\delta)})$. Then $T_{2}\in G_{[e_{1}]}$ and
$[\alpha T_{1}T_{2}]=[e_{1}+p^{r}e_{2}+a_{\nu+1}e_{\nu+1}+b_{\nu+2}e_{\nu+2}+p^{r_{2}}e_{\nu+3}+\tilde{\alpha}_{\delta}]$.
Since $r_{2}\geq r_{1}$ we have $r\leq r_{2}$. Then
$$
T_{3}=\left(
  \begin{array}{cccc|cccc|c}
   1&&&&0&&&&\\
   &1&&&&0&-p^{r_{2}-r}&&\\
   &&1&&&p^{r_{2}-r}&0&&\\
   &&&I^{(\nu-3)}&&&&&\\ \hline
   &&&&1&&&&\\
   &&&&&1&&&\\
   &&&&&&1&&\\
   &&&&&&&I^{(\nu-3)}&\\
   &&&&&&&&I^{(\delta)}\\
  \end{array}
\right)\in G_{[e_{1}]}
$$
and $[\alpha T_{1}T_{2}T_{3}]=[e_{1}+p^{r}e_{2}+a_{\nu+1}e_{\nu+1}+b_{\nu+2}e_{\nu+2}+\tilde{\alpha}_{\delta}]$,
where $b_{\nu+2}=0$ or $b_{\nu+2}=ax_{1}^{2}p^{t}$ for some $x_{1}\in R^{\ast}$, $r\leq t<s$ and $a\in\{1,z\}$.
Write $T_{0}=T_{1}T_{2}T_{3}$.

{\em Case\,1}: $\delta=0$. Then $[\alpha T_{0}]=[e_{1}+p^{r}e_{2}]$ if $b_{\nu+2}=0$; and $\mathrm{diag}(x_{1}I^{(\nu)},x_{1}^{-1}I^{(\nu)})\in G_{[e_{1}]}$ carries $[\alpha T_{0}]$ to $[e_{1}+p^{r}e_{2}-ap^{r+t}e_{\nu+1}+ap^{t}e_{\nu+2}]$ if $b_{\nu+2}=ax_{1}^{2}p^{t}$.

{\em Case\,2}: $\delta=1$. When $a_{2\nu+1}\ne 0$ we may write $a_{2\nu+1}=x_{2}p^{k}$, where $x_{2}\in R^{\ast}$ and $0<k<s$.

{\em Case\,2.1}: $a_{2\nu+1}=0$. Similar to {\em Case\,1}, there exists an element of $G_{[e_{1}]}$ carrying $[\alpha T_{0}]$ to $[e_{1}+p^{r}e_{2}]$ or $[e_{1}+p^{r}e_{2}-ap^{r+t}e_{\nu+1}+ap^{t}e_{\nu+2}]$.

{\em Case\,2.2}: $a_{2\nu+1}=x_{2}p^{k}$ and $r\leq k$. Then
$H_{2\nu+1,\Delta}(x_{2}p^{k-r}E_{2})\in G_{[e_{1}]}$
carries $[\alpha T_{0}]$ to $[e_{1}+p^{r}e_{2}+a_{\nu+1}e_{\nu+1}+(b_{\nu+2}+\frac{1}{2}\Delta x_{2}^{2}p^{2k-r})e_{\nu+2}]$, reduced to {\em Case\,2.1}.

{\em Case\,2.3}: $a_{2\nu+1}=x_{2}p^{k}$ and $r>k$. Let $c_{1}=\Delta x_{2}+p^{r-k}$, $c_{2}=\Delta x_{2}+p^{r-k}-\frac{1}{2}\Delta  ax_{1}^{2}p^{t-k}$ and $c_{3}=x_{2}-ax_{1}^{2}p^{t-k}$. Then $c_{1},c_{2},c_{3}\in R^{\ast}$ and
$T_{2\nu+1,\Delta}(E_{2})\in G_{[e_{1}]}$ carries $[\alpha T_{0}]$ to
$$
\begin{array}{ll}
[e_{1}+c_{1}p^{k}e_{2}+a_{\nu+1}e_{\nu+1}+x_{2}p^{k}e_{2\nu+1}] & \mbox{if $b_{\nu+2}=0$},\\[0.1cm]
[e_{1}+c_{2}p^{k}e_{2}+a_{\nu+1}e_{\nu+1}+b_{\nu+2}e_{\nu+2}+c_{3}p^{k}e_{2\nu+1}] & \mbox{if $b_{\nu+2}=ax_{1}^{2}p^{t}$}.\\
\end{array}
$$
If $b_{\nu+2}=0$, then
$\mathrm{diag}(1,c_{1}^{-1},I^{(\nu-1)},c_{1},I^{(\nu-1)})\in G_{[e_{1}]}$
carries $[\alpha T_{0} T_{2\nu+1,\Delta}(E_{2})]$ to $[e_{1}+p^{k}e_{2}+a_{\nu+1}e_{\nu+1}+x_{2}p^{k}e_{2\nu+1}]$, reduced to {\em Case\,2.2}.
If $b_{\nu+2}=ax_{1}^{2}p^{t}$, then $\mathrm{diag}(1,c_{2}^{-1},I^{(\nu-1)},c_{2},I^{(\nu-1)})\in G_{[e_{1}]}$ carries $[\alpha T_{0}T_{2\nu+1,\Delta}(E_{2})]$ to $[e_{1}+p^{k}e_{2}+a_{\nu+1}e_{\nu+1}+b_{\nu+2}c_{2}e_{\nu+2}+c_{3}p^{k}e_{2\nu+1}]$,
reduced to {\em Case\,2.2}.

{\em Case\,3}: $\delta=2$. When $a_{2\nu+i}\ne 0$ we write $a_{2\nu+i}=y_{i}p^{k_{i}}$, where $y_{i}\in R^{\ast}$, $0<k_{i}<s$ and $i=1,2$.

{\em Case\,3.1}: $a_{2\nu+1}=a_{2\nu+2}=0$. Similar to {\em Case\,1}, there exists an element of $G_{[e_{1}]}$ carrying $[\alpha T_{0}]$
to $[e_{1}+p^{r}e_{2}]$ or $[e_{1}+p^{r}e_{2}-ap^{r+t}e_{\nu+1}+ap^{t}e_{\nu+2}]$.

{\em Case\,3.2}: $a_{2\nu+1}\ne 0$ and $a_{2\nu+2}=0$. We may assume $r\leq k_{1}$. In fact, if $r>k_{1}$, then $d_{1}=y_{1}+p^{r-k_{1}}$, $d_{2}=y_{1}+p^{r-k_{1}}-\frac{1}{2}ax_{1}^{2}p^{t-k_{1}}$, $d_{3}=y_{1}-ax_{1}^{2}p^{t-k_{1}}$ are units, and
$T_{2\nu+2,\Delta}(E_{21})\in G_{[e_{1}]}$
carries $[\alpha T_{0}]$ to
$$
\begin{array}{ll}
[e_{1}+d_{1}p^{k_{1}}e_{2}+a_{\nu+1}e_{\nu+1}+y_{1}p^{k_{1}}e_{2\nu+1}] & \mbox{if $b_{\nu+2}=0$},\\[0.1cm]
[e_{1}+d_{2}p^{k_{1}}e_{2}+a_{\nu+1}e_{\nu+1}+b_{\nu+2}e_{\nu+2}+d_{3}p^{k_{1}}e_{2\nu+1}] & \mbox{if $b_{\nu+2}=ax_{1}^{2}p^{t}$}.\\
\end{array}
$$
If $b_{\nu+2}=0$, then
$\mathrm{diag}(1,d_{1}^{-1},I^{(\nu-1)},d_{1},I^{(\nu)})\in G_{[e_{1}]}$
carries $[\alpha T_{0}T_{2\nu+2,\Delta}(E_{21})]$ to $[e_{1}+p^{k_{1}}e_{2}+a_{\nu+1}e_{\nu+1}+y_{1}p^{k_{1}}e_{2\nu+1}]$, as desired.
If $b_{\nu+2}=ax_{1}^{2}p^{t}$, then $\mathrm{diag}(1,d_{2}^{-1},I^{(\nu-1)},d_{2},I^{(\nu)})\in G_{[e_{1}]}$ carries $[\alpha T_{0}T_{2\nu+2,\Delta}(E_{21})]$ to $$
[e_{1}+p^{k_{1}}e_{2}+a_{\nu+1}e_{\nu+1}+b_{\nu+2}d_{2}e_{\nu+2}+d_{3}p^{k_{1}}e_{2\nu+1}],
$$
as desired. Now assume $r\leq k_{1}$. Then
$H_{2\nu+2,\Delta}(y_{1}p^{k_{1}-r}E_{21})\in G_{[e_{1}]}$
carries $[\alpha T_{0}]$ to $[e_{1}+p^{r}e_{2}+a_{\nu+1}e_{\nu+1}+(b_{\nu+2}+\frac{1}{2}y_{1}^{2}p^{2k_{1}-r})e_{\nu+2}]$, reduced to {\em Case\,3.1}.

{\em Case\,3.3}: $a_{2\nu+1}=0$ and $a_{2\nu+2}\ne 0$. Similar to {\em Case\,3.2}, there exists an element of $G_{[e_{1}]}$ carrying $[\alpha T_{0}]$
to $[e_{1}+p^{r}e_{2}]$ or $[e_{1}+p^{r}e_{2}-ap^{r+t}e_{\nu+1}+ap^{t}e_{\nu+2}]$.

{\em Case\,3.4}: $a_{2\nu+1}\ne 0$ and $a_{2\nu+2}\ne 0$.

{\em Case\,3.4.1}: $\min\{r,k_{1},k_{2}\}=r$. Then
$H_{2\nu+2,\Delta}(y_{1}p^{k_{1}-r}E_{21})\in G_{[e_{1}]}$
carries $[\alpha T_{0}]$ to $[e_{1}+p^{r}e_{2}+a_{\nu+1}e_{\nu+1}+(b_{\nu+2}+\frac{1}{2}y_{1}^{2}p^{2k_{1}-r})e_{\nu+2}+a_{2\nu+2}e_{2\nu+2}]$, reduced to {\em Case\,3.3}.

{\em Case\,3.4.2}: $\min\{r,k_{1},k_{2}\}=k_{1}<r$. Then $d_{1}=y_{1}+p^{r-k_{1}}$, $d_{2}=y_{1}+p^{r-k_{1}}-\frac{1}{2}ax_{1}^{2}p^{t-k_{1}}$, $d_{3}=y_{1}-ax_{1}^{2}p^{t-k_{1}}$ are units, and
$T_{2\nu+2,\Delta}(E_{21})\in G_{[e_{1}]}$
carries $[\alpha T_{0}]$ to
$$
\begin{array}{ll}
[e_{1}+d_{1}p^{k_{1}}e_{2}+a_{\nu+1}e_{\nu+1}+y_{1}p^{k_{1}}e_{2\nu+1}+y_{2}p^{k_{2}}e_{2\nu+2}] & \mbox{if $b_{\nu+2}=0$},\\[0.1cm]
[e_{1}+d_{2}p^{k_{1}}e_{2}+a_{\nu+1}e_{\nu+1}+b_{\nu+2}e_{\nu+2}+d_{3}p^{k_{1}}e_{2\nu+1}+y_{2}p^{k_{2}}e_{2\nu+2}] & \mbox{if $b_{\nu+2}=ax_{1}^{2}p^{t}$}.\\
\end{array}
$$
If $b_{\nu+2}=0$, then
$\mathrm{diag}(1,d_{1}^{-1},I^{(\nu-1)},d_{1},I^{(\nu)})\in G_{[e_{1}]}$
carries $[\alpha T_{0}T_{2\nu+2,\Delta}(E_{21})]$ to
$[e_{1}+p^{k_{1}}e_{2}+a_{\nu+1}e_{\nu+1}+y_{1}p^{k_{1}}e_{2\nu+1}+y_{2}p^{k_{2}}e_{2\nu+2}]$,
reduced to {\em Case\,3.4.1}.
If $b_{\nu+2}=ax_{1}^{2}p^{t}$, then
$\mathrm{diag}(1,d_{2}^{-1}$, $I^{(\nu-1)},d_{2},I^{(\nu)})\in G_{[e_{1}]}$
carries $[\alpha T_{0}T_{2\nu+2,\Delta}(E_{21})]$ to
$[e_{1}+p^{k_{1}}e_{2}+a_{\nu+1}e_{\nu+1}+b_{\nu+2}d_{2}e_{\nu+2}+d_{3}p^{k_{1}}e_{2\nu+1}+y_{2}p^{k_{2}}e_{2\nu+2}]$,
reduced to {\em Case\,3.4.1}.

{\em Case\,3.4.3}: $\min\{r,k_{1},k_{2}\}=k_{2}<r$. Similar to {\em Case\,3.4.2}, there exists an element of $G_{[e_{1}]}$ carrying $[\alpha T_{0}]$
to $[e_{1}+p^{r}e_{2}]$ or $[e_{1}+p^{r}e_{2}-ap^{r+t}e_{\nu+1}+ap^{t}e_{\nu+2}]$.
$\qed$

Note that $\{[e_{1}]\}$ is the trivial orbit of $G_{[e_{1}]}$ on $\Gamma^{2\nu+\delta}$.
Combining Lemmas~\ref{lem3.2} and \ref{lem3.3}, we obtain the following result.

\begin{thm}\label{thm3.4}
The nontrivial orbits of $G_{[e_{1}]}$ on $\Gamma^{2+\delta}$ have the following representatives:
\begin{equation}
\begin{array}{lll}\label{formula1.x}
&[e_{2}]& \mbox{if $\delta=0$},\\[0.1cm]
&[e_{2}], [1,-\frac{1}{2}\Delta p^{2r},p^{r}]& \mbox{if $\delta=1$},\\[0.1cm]
&[e_{2}], [1,-\frac{1}{2}p^{2r},p^{r},0], [1,\frac{1}{2}zp^{2r},0,p^{r}], [1,-\frac{1}{2}(b^{2}-z)p^{2r},bp^{r},p^{r}]& \mbox{if $\delta=2$},\\
\end{array}
\end{equation}
where $b\in R^{\ast}$, $0<r<s$ and $z$ is a fixed non-square element of $R^{\ast}$.
When $\nu\geq 2$, the nontrivial orbits of $G_{[e_{1}]}$ on $\Gamma^{2\nu+\delta}$ have the following representatives:
\begin{equation}\label{formula2.x}
[e_{\nu+1}],[e_{2}],[e_{2}+p^{r}e_{\nu+1}],[e_{1}+p^{r}e_{2}],[e_{1}+p^{r}e_{2}-ap^{r+t}e_{\nu+1}+ap^{t}e_{\nu+2}]
\end{equation}
where  $1\leq r\leq t\leq s-1$ and $a\in\{1,z\}$.
$\qed$
\end{thm}

It should be noted that two distinct representatives listed in Theorem~\ref{thm3.4} may fall into the same orbit of $G_{[e_{1}]}$. For example, when $(\nu,\delta)=(1,2)$ and $-1\notin {R^{\ast}}^{2}$, pick $x\in R^{\ast}$ such that $zx^2=-1$. Then
$$
\left(
  \begin{array}{cc|cc}
   x&&&\\
   &x^{-1}&&\\ \hline
   &&0&x\\
   &&xz&0\\
  \end{array}
\right)\in G_{[e_{1}]}
$$
carries $[1,-\frac{1}{2}p^{2r},p^{r},0]$ to $[1,\frac{1}{2}zp^{2r},0,p^{r}]$.

\section{Quasi-strongly regular graphs}

A {\em strongly regular graph} with parameters $(n,k,\lambda,\mu)$ is an undirected regular graph with valency $k$ on $n$ vertices
such that each pair of adjacent vertices has $\lambda$ common neighbors, and each pair of nonadjacent
vertices has $\mu$ common neighbors.
As a generalization of strongly regular graphs, quasi-strongly regular graphs were discussed by \cite{Gold} and \cite{gol1}. Let $c_1,c_2,\dots,c_d$ be distinct non-negative integers. A connected graph of valency $k$ on $n$ vertices is
{\em quasi-strongly regular} of grade $d$ with parameters $(n,k,\lambda;c_1,c_2,\dots,c_d)$ if any two adjacent vertices
have $\lambda$ common neighbors, and any two non-adjacent vertices have $c_i$ common neighbors
for some $i$.
A quasi-strongly regular graph with grade $1$ is strongly regular. As pointed in \cite{Gold}, quasi-strongly regular graphs with grade $2$ are important.
Since $\Gamma^{2\nu+\delta}$ is an arc transitive graph, it is quasi-strongly regular.
In this section we shall show that $\Gamma^{2\nu+\delta}$ has grade at most two, and compute all its parameters.

\begin{prop}~\label{prop}
The numbers of $(x,y)\in R\times R$ with $x^{2}-zy^{2}\notin R^{\ast}$ is $p^{2m(s-1)}$, where $z$ is a fixed non-square element of $R^{\ast}$.
\end{prop}
\proof
Let $\Omega=\{(x,y)\in R\times R\mid x^{2}-zy^{2}\notin R^{\ast}\}$. For $r\in R$, write $\bar r=r+pR$.
Assume that $(x,y)\in\Omega$. Then $\bar{x}^{2}-\bar{z}\bar{y}^{2}=\bar{0}$ in the quotient field
$R/pR$. By \cite[Chapter 6]{Wan}, the matrix $\mathrm{diag}(\bar{1},-\bar{z})$ over $R/pR$ is definite. So $\bar{x}=\bar{y}=\bar{0}$, i.e., $x,y\in pR$. It follows that $\Omega=\{(x,y)\mid x,y\in pR\}$ and $|\Omega|=p^{2m(s-1)}$.
$\qed$

For convenience we write $\tilde{x}_{\delta}=\emptyset (\mathrm{disappear})$, $-\frac{1}{2}\Delta x_{2\nu+1}^{2}$ or $-\frac{1}{2}(x_{2\nu+1}^{2}-zx_{2\nu+2}^{2})$ according to $\delta=0,1$ or 2, respectively.

\begin{lemma}\label{lem4.A}
$\Gamma^{2\nu+\delta}$ is a regular graph of valency $p^{ms(2\nu-2+\delta)}$ on
$$
\frac{(p^{m\nu}-1)(p^{m(\nu+\delta-1)}+1)p^{m(s-1)(2\nu-2+\delta)}}{p^{m}-1}
$$
many vertices.
\end{lemma}
\proof
Each vertex adjacent to $[e_{1}]$ is of the form $[x_{1},\ldots,x_{\nu},1,x_{\nu+2},\ldots,x_{2\nu+\delta}]$, where
$x_{1}=-x_{2}x_{\nu+2}-\cdots-x_{\nu}x_{2\nu}+\tilde{x}_{\delta}$. By Theorem~\ref{thm3.1}, the valency of $\Gamma^{2\nu+\delta}$ is $p^{ms(2\nu-2+\delta)}$.

Now we compute the number of vertices of $\Gamma^{2\nu+\delta}$. Note that there exists a unit among $x_{1},x_{2},\ldots,x_{2\nu}$ for any vertex $[x_{1},x_{2},\ldots,x_{2\nu+\delta}]$ of $\Gamma^{2\nu+\delta}$. For $1\leq i\leq 2\nu$, let
$$
\Omega_{i}=\{[x_{1},\ldots,x_{i-1},1,x_{i+1},\ldots,x_{2\nu+\delta}]\in V(\Gamma^{2\nu+\delta}) \mid x_{1},\ldots,x_{i-1}\notin R^{\ast}\}.
$$
Then $\{\Omega_{1},\Omega_{2},\ldots,\Omega_{2\nu}\}$ is a partition of $\tilde{R}^{2\nu+\delta}$.
Any element $[x_{1},\ldots,x_{i-1},1,x_{i+1},\ldots,x_{2\nu+\delta}]$ of $\Omega_{i}$ satisfies
$$
x_{\nu+i}=-x_{1}x_{\nu+1}-\cdots-x_{i-1}x_{\nu+i-1}-x_{i+1}x_{\nu+i+1}-\cdots-x_{\nu}x_{2\nu}+\tilde{x}_{\delta},
$$
where $1\leq i\leq \nu$ and $x_{1},\ldots,x_{i-1}\notin R^{\ast}$. Since there are $p^{m(s-1)(i-1)}$ choices for $(x_{1},\ldots,x_{i-1})$; and for a given $(x_1,\ldots,x_{i-1})$, there are $p^{ms(2\nu+\delta-i-1)}$ choices for
$$
(x_{i+1},\ldots,x_{\nu+i-1},x_{\nu+i+1},\ldots,x_{2\nu+\delta}),
$$
we have
$|\Omega_{i}|=p^{m(s(2\nu-2+\delta)+1-i)}$.
Any element $[x_{1},\ldots,x_{\nu+i-1},1,x_{\nu+i+1},\ldots,x_{2\nu+\delta}]$ of $\Omega_{\nu+i}$ satisfies
$$
x_{i}=-x_{1}x_{\nu+1}-\cdots-x_{i-1}x_{\nu+i-1}-x_{i+1}x_{\nu+i+1}-\cdots-x_{\nu}x_{2\nu}+\tilde{x}_{\delta},
$$
where $1\leq i\leq \nu$ and $x_{1},\ldots,x_{\nu+i-1}\notin R^{\ast}$. Note that $x_{2\nu+1}\notin R^{\ast}$ if $\delta=1$, and $x_{2\nu+1},x_{2\nu+2}\notin R^{\ast}$ if $\delta=2$.
There are $p^{ms(\nu-i)}$ choices for $(x_{\nu+i+1},\ldots,x_{2\nu})$ and there are $p^{m(s-1)(\nu+i-2+\delta)}$ choices for
$(x_{1},\ldots,x_{i-1},x_{i+1},\ldots,x_{\nu+i-1},\tilde{x})$, where $\tilde{x}=\emptyset,x_{2\nu+1}$ or $(x_{2\nu+1},x_{2\nu+2})$ according to
$\delta=0,1$ or 2, respectively. We have $|\Omega_{\nu+i}|=p^{m(s\nu+(s-1)(\nu-2+\delta)-i)}$.
Therefore, the desired result follows.
$\qed$

\begin{thm}\label{thm4.3}
$\Gamma^{2+\delta}$ is strongly regular with parameters $(n,k,\lambda,\mu)$, where $n,k$ are given by Lemma~\ref{lem4.A},
$\lambda=(p^{m\delta}-1)p^{m\delta (s-1)}$ and $\mu=\lceil\delta/2\rceil p^{ms\delta}$.
\end{thm}
\proof
Let $[\alpha]$ and  $[\beta]$ be any two adjacent vertices of $\Gamma^{2+\delta}$.
By Theorems~\ref{thm3.1} and \ref{thm3.4} we may assume that $[\alpha]=[e_1]$ and $[\beta]=[e_{2}]$.
If $\delta=0$, then $[\alpha]$ and $[\beta]$ have no common neighbors.
If $\delta=1$, then any common neighbor of $[\alpha]$ and $[\beta]$ is of the form $[1,-\frac{1}{2}\Delta x^{2},x]$, where $x\in R^{\ast}$.
So $\lambda=(p^{m}-1)p^{m(s-1)}$.
If $\delta=2$, then any common neighbor of $[\alpha]$ and $[\beta]$ is of the form $[1,-\frac{1}{2}(x^{2}-zy^{2}),x,y]$, where $x^{2}-zy^{2}\in R^{\ast}$. By Proposition~\ref{prop}, $\lambda=(p^{2m}-1)p^{2m(s-1)}$.

Let $[\alpha]$ and  $[\gamma]$ be any two non-adjacent vertices of $\Gamma^{2+\delta}$.
Note that $[\alpha]$ and $[\gamma]$ have no common neighbors when $\delta=0$.
When $\delta\geq 1$, by Theorems~\ref{thm3.1} and \ref{thm3.4}, we may choose $[\alpha]=[e_{1}]$ and $[\gamma]$ to be one of the vertices listed in (\ref{formula1.x}) except $[e_{2}]$.
If $\delta=1$, then any common neighbor of $[\alpha]$ and $[\gamma]$ is of the form $[-\frac{1}{2}\Delta x_{3}^{2},1,x_{3}]$. So $\mu=p^{ms}$.
Assume that $\delta=2$. If $[\gamma]=[1,-\frac{1}{2}p^{2r},p^{r},0]$, then any common neighbor of $[\alpha]$ and $[\gamma]$ is of the form
$[-\frac{1}{2}(x_{3}^{2}-zx_{4}^{2}),1,x_{3},x_{4}]$. So $[\alpha]$ and $[\gamma]$ have $p^{2ms}$ common neighbors.
Similarly, if $[\gamma]=[1,\frac{1}{2}zp^{2r},0,p^{r}]$ or $[1,-\frac{1}{2}(b^{2}-z)p^{2r},bp^{r},p^{r}]$, then $[\alpha]$ and $[\gamma]$ have $p^{2ms}$ common neighbors.
$\qed$

\begin{thm}\label{thm4.4}
If $\nu\geq 2$, then $\Gamma^{2\nu+\delta}$ is quasi-strongly regular with parameters $(n,k,\lambda;c_1,c_2)$,
where $n,k$ are given by Lemma~\ref{lem4.A},
$\lambda=(p^{m}-1)(1+(\delta-1)p^{m(1-\nu-\frac{\delta}{2})})p^{m(2s\nu-2s-1+s\delta)}$, $c_{1}=(p^{m}-1)p^{m(s(2\nu-2+\delta)-1)}$ and $c_{2}=p^{ms(2\nu-2+\delta)}$.
\end{thm}
\proof
First we compute $\lambda$. Let $[\alpha]$ and $[\beta]$ be any two adjacent vertices of $\Gamma^{2\nu+\delta}$.
By Theorems~\ref{thm3.1} and \ref{thm3.4} we may assume that $[\alpha]=[e_{1}]$ and $[\beta]=[e_{\nu+1}]$.
Then any common neighbor of $[\alpha]$ and $[\beta]$ is of the form $[1,x_{2},\ldots,x_{2\nu+\delta}]$, where $x_{\nu+1}\in R^{\ast}$ and
$x_{\nu+1}=-x_{2}x_{\nu+2}-\cdots-x_{\nu}x_{2\nu}+\tilde{x}_{\delta}$.  Note that there exists a unit among $x_{2},\ldots,x_{\nu},\tilde{x}_{\delta}$.
Let
$$
\Omega_{i}=\left\{\begin{array}{ll}
\{[1,x_{2},\ldots,x_{2\nu+\delta}]\in V(\Gamma^{2\nu+\delta}) \mid x_{i}\in R^{\ast},x_{2},\ldots,x_{i-1}\notin R^{\ast}\}& \mbox{$2\leq i\leq \nu$},\\[0.1cm]
\{[1,x_{2},\ldots,x_{2\nu+\delta}]\in V(\Gamma^{2\nu+\delta}) \mid \tilde{x}_{\delta}\in R^{\ast},x_{2},\ldots,x_{\nu}\notin R^{\ast}\}& \mbox{$i=2\nu+\delta$}.\\[0.1cm]
\end{array}\right.
$$
For $2\leq i\leq \nu$, any element $[1,x_{2},\ldots,x_{2\nu+\delta}]$ of $\Omega_{i}$ satisfies
$$
x_{\nu+i}=-x_{i}^{-1}(x_{\nu+1}+x_{2}x_{\nu+2}+\cdots+x_{i-1}x_{\nu+i-1}+x_{i+1}x_{\nu+i+1}+\cdots+x_{\nu}x_{2\nu}-\tilde{x}_{\delta}).
$$
Since there are $(p^{m}-1)^{2}p^{2m(s-1)}$ choices for $(x_{i},x_{\nu+1})$, there are $p^{m(s-1)(i-2)}$ choices for $(x_{2},\ldots,x_{i-1})$, and for a given $(x_{2},\ldots,x_{i-1})$ there are  $p^{ms(2\nu-2+\delta-i)}$ choices for
$$
(x_{i+1},\ldots,x_{\nu},x_{\nu+2},\ldots,x_{\nu+i-1},x_{\nu+i+1},\ldots,x_{2\nu+\delta}),
$$
the size of $\Omega_i$ is $(p^{m}-1)^{2}p^{m(s(2\nu-2+\delta)-i)}$.
Any element $[1,x_{2},\ldots,x_{2\nu+\delta}]$ of $\Omega_{2\nu+\delta}$ satisfies
$x_{\nu+1}=-x_{2}x_{\nu+2}-\cdots-x_{\nu}x_{2\nu}+\tilde{x}_{\delta}$.
Since there are $(p^{\delta m}-1)p^{\delta m(s-1)}$ choices for $(x_{2\nu+1},\ldots,x_{2\nu+\delta})$ by Proposition~\ref{prop},
there are $p^{m(s-1)(\nu-1)}$ choices for $(x_{2},\ldots,x_{\nu})$, and  there are  $p^{ms(\nu-1)}$ choices for $(x_{\nu+2},\ldots,x_{2\nu})$,
the size of $\Omega_{2\nu+\delta}$ is $(p^{\delta m}-1)p^{m((s-1)(\nu+\delta-1)+s(\nu-1))}$.
Note that $\lambda=|\Omega_{2}|+\cdots+|\Omega_{\nu}|+\vartheta$, where $\vartheta=0$, $|\Omega_{2\nu+1}|$ or $|\Omega_{2\nu+2}|$ according to $\delta=0$, 1 or 2, respectively. The desired $\lambda$ follows.

Next we compute $c_{1}$ and $c_{2}$. Let $[\alpha]$ and $[\gamma]$ be any two non-adjacent vertices of $\Gamma^{2\nu+\delta}$.
By Theorems~\ref{thm3.1} and \ref{thm3.4}, we may choose $[\alpha]=[e_{1}]$ and $[\gamma]$ to be one of the vertices listed in (\ref{formula2.x})
except $[e_{\nu+1}]$.
If $[\gamma]=[e_{2}]$, then any common neighbor of $[\alpha]$ and $[\gamma]$ is of the form
$[x_{1},\ldots,x_{\nu},1,x_{\nu+2},\ldots,x_{2\nu+\delta}]$, where $x_{\nu+2}\in R^{\ast}$ and
$x_{2}=-x_{\nu+2}^{-1}(x_{1}+x_{3}x_{\nu+3}+\cdots+x_{\nu}x_{2\nu}-\tilde{x}_{\delta})$.
So the number of such vertices is $(p^{m}-1)p^{m(s(2\nu-2+\delta)-1)}$.
Similarly, if $[\gamma]=[e_{2}+p^{r}e_{\nu+1}]$,
then $[\alpha]$ and $[\gamma]$ have $(p^{m}-1)p^{m(s(2\nu-2+\delta)-1)}$ common neighbors.
If $[\gamma]=[e_{1}+p^{r}e_{2}]$, then any common neighbor of $[\alpha]$ and $[\gamma]$ is of the form
$[x_{1},\ldots,x_{\nu},1,x_{\nu+2},\ldots,x_{2\nu+\delta}]$, where $1+x_{\nu+2}p^{r}\in R^{\ast}$ and
$x_{1}+x_{2}x_{\nu+2}+x_{3}x_{\nu+3}+\cdots+x_{\nu}x_{2\nu}-\tilde{x}_{\delta}=0$.
So any common neighbor $[x_{1},\ldots,x_{\nu},1,x_{\nu+2},\ldots,x_{2\nu+\delta}]$ of $[\alpha]$ and $[\gamma]$ satisfies
$x_{1}=-(x_{2}x_{\nu+2}+x_{3}x_{\nu+3}+\cdots+x_{\nu}x_{2\nu}-\tilde{x}_{\delta})$.
It follows that the number of such vertices is $p^{ms(2\nu-2+\delta)}$.
Similarly, if $[\gamma]=[e_{1}+p^{r}e_{2}-ap^{r+t}e_{\nu+1}+ap^{t}e_{\nu+2}]$, then
$[\alpha]$ and $[\gamma]$ have $p^{ms(2\nu-2+\delta)}$ common neighbors.
$\qed$

Theorem~\ref{thm4.3} and \ref{thm4.4} show that the diameter of $\Gamma^{2\nu+\delta}$ is 2 unless $(\nu,\delta)=(1,0)$. In the case of $(\nu,\delta)=(1,0)$, it is a clique with two vertices.

As another generalization of strongly regular graphs, Erickson et al. \cite{Erickson} introduced Deza graphs, which were firstly
introduced in a slightly more restricted form by Deza and Deza \cite{DD}. A regular graph with valency $k$ on $n$ vertices is said
to be a $(n,k,b,c)$-{\it Deza graph} if any two distinct vertices $x$ and $y$ have $b$ or $c$ common adjacent vertices. A Deza graph
with diameter two is said to be a {\it strictly Deza graph} if it is not strongly regular.
For $\nu\geq 2$, Theorem~\ref{thm4.4} shows that $\Gamma^{2\nu+1}$ is a strictly Deza graph with parameters
$$
((p^{2m\nu}-1)p^{m(s-1)(2\nu-1)}/(p^{m}-1),\;p^{ms(2\nu-1)},\;(p^{m}-1)p^{m(2s\nu-s-1)},\;p^{ms(2\nu-1)}).
$$

\section*{Acknowledgment}
This research is supported by the Scientific Research Fund of Hunan Provincial Education Department (12A058),
NSFC (11271047), NSF of Hebei Province (A2012408003, A2013408009), TPF-2011-11 of Hebei
Province, NSF of Hebei Education Department (ZH2012082) and the Fundamental Research Funds for the Central University of China.

\end{CJK*}

\end{document}